\documentclass[10pt]{article}

\usepackage{amsmath,amsthm,amsfonts,amssymb,amscd}
\newtheorem {Lemma}{Lemma}
\newtheorem {Theorem}{Theorem}
\newenvironment {Proof} {\noindent {\bf Proof.}}{\quad $\square$\par\vspace{3mm}}
\allowdisplaybreaks [4]

\begin{document}

\title{\bf A note on skew spectrum of graphs}

\author{Yanna Wang,  Bo Zhou\footnote{Corresponding author. E-mail address: zhoubo@scnu.edu.cn}
\\[1mm]
Department of Mathematics, South China Normal University,\\
Guangzhou 510631, P.~R. China
}

\date{}

\maketitle

\begin{abstract}
We give some properties of skew spectrum of a graph, especially, we
answer  negatively a problem  concerning the skew characteristic
polynomial and matching polynomial in [M. Cavers et al.,
Skew-adjacency matrices of graphs, Linear Algebra Appl. 436 (2012)
 4512--4529].
\end{abstract}

\section{Introduction}

We consider simple graphs. Let $G$ be a graph with vertex set $V(G)$
and edge set $E(G)$. An orientation of  $G$ is a sign-valued
function $\sigma$ on the set of ordered pairs $\{(i, j), (j,
i)|ij\in E(G)\}$ that specifies an orientation  to each edge $ij$ of
$G$: If $ij\in E(G)$, then we take  $\sigma(i, j)=1$   when
$i\rightarrow j$ and  $\sigma(i,j)=-1$ when $j\rightarrow i$. The
resulting oriented graph is denoted by $G^{\sigma}$. Both $\sigma$
and $G^{\sigma}$ are called orientations of $G$.

The skew adjacency matrix $S^{\sigma}=S(G^{\sigma})$ of $G^{\sigma}$
is the $\{0,1,-1\}$-matrix with $(i, j)$-entry equal to $\sigma(i,
j)$ if $ij\in E(G^{\sigma})$ and $0$ otherwise. If there is no
confusion, we simply write $S = [s_{i,j}]$ for $S^{\sigma}$. Thus
$s_{i, j}=1$ if $ij\in E(G^{\sigma})$, $-1$ if $ji\in
E(G^{\sigma})$, and $0$ otherwise. The (skew) characteristic
polynomial of $S=S^{\sigma}$ is
\[
p_S(x)=\det(xI-S)=x^{n}+s_{1}x^{n-1}+\cdots+s_{n-1}x+s_{n},
\]
where $n=|V(G)|$. Let $\rho(D)$ be the spectral radius of a square
matrix $D$, i.e., the largest modulus of the eigenvalues of $D$. The
spectral radius of $G$ is the spectral radius of its adjacency
matrix. The maximum skew spectral radius of $G$ is defined as
$\rho_{s}(G) = \max_{S} \rho(S)$, where the maximum is taken over
all of the skew adjacency matrices $S$ of $G$.

An odd-cycle graph is a graph with no  even cycles (cycles of even
lengths). In particular, a tree is an odd-cycle graph.

Let $G$ be a graph with $n$ vertices. Let $m_{k}(G)$ be the number
of matchings in $G$ that cover $k$ vertices. Obviously, $m_k(G)=0$
if $k$ is odd. The matching polynomial of $G$ is defined as
\[
m(G,x)=\sum_{k=0}^{n}(-1)^{\frac{k}{2}}m_{k}(G)x^{n-k},
\]
where $m_{0}(G)=1$.

Let $G$ be a graph on $n$ vertices. After showing that $G$ is an
odd-cycle graph if and only if $p_S(x)=(-i)^{n}m(G,ix)$ for all skew
adjacency matrices $S$ of $G$ (see \cite[Lemma~5.4]{Ca}), Cavers et
al. \cite{Ca} posed the following question:

\noindent {\bf Problem 1.} If $p_{S}(x)=(-i)^{n}m(G,ix)$ for some
skew adjacency matrix $S$ of $G$, must $G$ be an odd-cycle graph?

After showing that if $G$ is an odd-cycle graph, then
$\rho_s(G)=\rho(S)$ for all skew adjacency matrices $S$ of $G$
(\cite[Lemma~6.2]{Ca}), Cavers et al. \cite{Ca} posed the following
question:

\noindent {\bf Problem 2.} If $G$ is a connected graph and $\rho(S)$
is the same for all skew adjacency matrices $S$ of $G$, must $G$ be
an odd-cycle graph?

In this note we answer Problem 1 negatively by constructing a class
of graphs, and when $G$ is a connected bipartite graph we answer
Problem 2 affirmatively.

\section{Preliminaries}

Let $\mathcal{U}_{k}$ be the set of all collections $U$ of (vertex)
 disjoint edges and even cycles in $G$ that cover $k$ vertices ($\mathcal{U}_{k}^e$ was used for this set in \cite{Ca}).  A routing $\overrightarrow{U}$ of $U\in \mathcal{U}_{k}$
is obtained by replacing each edge in $U$ by a digon and each (even)
cycle in $U$ by a dicycle. For  an orientation $\sigma$ of a graph
$G$ and a routing $\overrightarrow{U}$  of $U\in \mathcal{U}_{k}$,
let $\sigma (\overrightarrow{U})= \prod_{(i, j)\in
E(\overrightarrow{U})}\sigma(i, j)$. We say that
$\overrightarrow{U}$ is positively  (resp. negatively) oriented
relative to $\sigma$ if $\sigma (\overrightarrow{U})=1$ (resp.
$\sigma (\overrightarrow{U})=-1$). For $U\in \mathcal{U}_{k}$, let
$c^{+}(U)$ (resp. $c^{-}(U)$) be the number of  cycles in $U$ that
are positively (resp. negatively ) oriented relative to $\sigma$
when $U$ is given a routing $\overrightarrow{U}$. Then
$c(U)=c^+(U)+c^-(U)$ is the (total) number of cycles of $U$.

\begin{Lemma} \label{lm2.1} \textnormal{\cite[eq.~(8)]{Ca}}
Let $S$ be a  skew adjacency matrix of $G$.
Then $s_{k}=0$
if $k$ is odd and
\[
s_{k}=m_{k}(G)+\sum_{U\in \mathcal{U}_{k}\atop
c(U)>0}(-1)^{c^{+}(U)}2^{c(U)} \mbox{ if $k$ is even}.
\]
\end{Lemma}


The following lemma is obtained from
parts 2 and 3 of Lemma 6.2 in \cite{Ca}.

\begin{Lemma} \label{lm2.5} Let $G$ be a connected bipartite graph,
\[
A=\left[
\begin{array}{ccc}
 O & B  \\
 B^{\top} & O
\end{array}
\right]
\]
the adjacency matrix of $G$, and
\[
S=\left[
\begin{array}{ccc}
 O & B  \\
 -B^{\top} & O
\end{array}
\right]
\]
and
\[
\tilde{S}=\left[
\begin{array}{ccc}
 O & \widetilde{B}  \\
 -\widetilde{B}^{\top} & O
\end{array}
\right]
\]
two skew adjacency matrices of $G$. Then $\rho(A)=\rho_{s}(G)$, and
$\rho(S)=\rho(\widetilde{S})$ if and only if $\widetilde{S}=DSD^{-1}$ for
some $\{-1, 1\}$-diagonal matrix $D$.
\end{Lemma}


\begin{Lemma} \label{lm2.2} \textnormal{\cite[Theorem 4.2]{Ca}}
The skew adjacency matrices of a graph $G$ are all cospectral if and only if $G$ is an odd-cycle graph.
\end{Lemma}

%
%
%

\section{Results}

First we give a negative answer to Problem 1.

\begin{Theorem} \label{Th1} For integer $m\ge 2$, let $G$ be the graph consisting
of two $2m$-vertex cycles $C_1$ and $C_2$ with exactly one common
vertex. Let $\sigma$ be an orientation of $G$ such that the cycle
$C_{1}$ $($resp. ${C_{2}})$ is positively $($resp. negatively$)$
oriented relative to $\sigma$. Let $S=S(G^{\sigma})$ and let
$n=4m-1$. Then $p_{S}(x)=(-i)^{n}m(G,ix)$.

\end{Theorem}

\begin{Proof}
It is sufficient to show that $s_{k}=m_{k}(G)$ for even $k$. By
Lemma~\ref{lm2.1}, we only need to show that
\[
\sum_{U\in \mathcal{U}_{k}\atop c(U)>0}(-1)^{c^{+}(U)}2^{c(U)}=0
\mbox{ for even $k$}.
\]
This is obvious when $k<2m$. Suppose that $k$ is even with  $2m\leq
k\leq4m-2$. Let $C_1=v_1v_2\dots v_{2m}v_1$ and
$C_2=v_{1}'v_{2}'\dots v_{2m}' v_1'$ with $v_1=v_1'$.

Let $\mathcal{U}_{k}^1$ be the subset of $\mathcal{U}_{k}$
consisting of $C_1$ and $\frac{1}{2}(k-2m)$ disjoint edges in $C_2$
and $\mathcal{U}_{k}^2$ the subset of $\mathcal{U}_{k}$ consisting
of $C_2$ and $\frac{1}{2}(k-2m)$ disjoint edges in $C_1$. Obviously, $\mathcal{U}_{k}^1\cap \mathcal{U}_{k}^2=\emptyset$.  For any
$U\in \mathcal{U}_{k}$ with $c(U)>0$, $U\in \mathcal{U}_{k}^1$ or
$U\in \mathcal{U}_{k}^2$. There is a bijection  from
$\mathcal{U}_{k}^1$ to $\mathcal{U}_{k}^2$ which maps $U\in
\mathcal{U}_{k}^1$ consisting of $C_1$ and $\frac{1}{2}(k-2m)$
disjoint edges, say $v_{i_1}'v_{i_1+1}',\dots,
v_{i_s}'v_{i_s+1}'$ in $C_2$ to $U'\in \mathcal{U}_{k}^2$ consisting of $C_2$
and $\frac{1}{2}(k-2m)$ disjoint edges $v_{i_1}v_{i_1+1},\dots,
v_{i_s}v_{i_s+1}$  in $C_1$, where $s=\frac{1}{2}(k-2m)$ and $2\le
i_1<\dots<i_s\le 2m-1$. Thus
$|\mathcal{U}_{k}^1|=|\mathcal{U}_{k}^2|$. Note that
\[
\sum_{U\in
\mathcal{U}_{k}^1}(-1)^{c^{+}(U)}2^{c(U)}=(-1)^{1}\cdot2^{1}\cdot
|\mathcal{U}_{k}^1|
\]
and
\[
\sum_{U\in
\mathcal{U}_{k}^2}(-1)^{c^{+}(U)}2^{c(U)}=(-1)^{0}\cdot2^{1}\cdot
|\mathcal{U}_{k}^2|.
\]
Thus \[ \sum_{U\in \mathcal{U}_{k}\atop
c(U)>0}(-1)^{c^{+}(U)}2^{c(U)}=\sum_{U\in
\mathcal{U}_{k}^1}(-1)^{c^{+}(U)}2^{c(U)}+\sum_{U\in
\mathcal{U}_{k}^2}(-1)^{c^{+}(U)}2^{c(U)}=0,
\]
as desired.
\end{Proof}

See Fig.~1 for an example with $7$ vertices for Theorem \ref{Th1}
and its proof.

\setlength{\unitlength} {0.7mm}\label{fig.1}

\begin{center}

\begin{picture}(100,80)

\put(50,47){\circle* {2.5}}         \put(50,50) {$v_{1}(v_{1}^{'})$}

\put(50,47){\line(1,0) {25}} \put(50,47){\vector(1,0) {15}}

\put(75,47){\circle* {2.5}}  \put(75,50) {$v_{2}^{'}$}

\put(50,22){\line(0,1) {25}} \put(50,22){\vector(0,1) {15}}

\put(50,22){\circle* {2.5}} \put(48,15) {$v_{4}^{'}$}

\put(75,47){\line(0,-1) {25}} \put(75,47){\vector(0,-1) {15}}

\put(75,22){\circle* {2.5}}  \put(75,15) {$v_{3}^{'}$}

\put(50,22){\line(1,0) {25}}\put(50,22){\vector(1,0) {15}}

\put(50,47){\line(-1,0) {25}} \put(50,47){\vector(-1,0) {15}}
\put(25,47){\circle* {2.5}} \put(25,42) {$v_{2}$}

\put(50,72){\line(0,-1) {25}}  \put(50,72){\vector(0,-1) {15}}
\put(50,72){\circle* {2.5}} \put(50,74) {$v_{4}$}

\put(25,72){\line(1,0) {25}} \put(25,72){\vector(1,0) {15}}
\put(25,72){\circle* {2.5}} \put(25,74) {$v_{3}$}

\put(25,47){\line(0,1) {25}} \put(25,47){\vector(0,1) {15}}

\end{picture}
\vskip -7mm
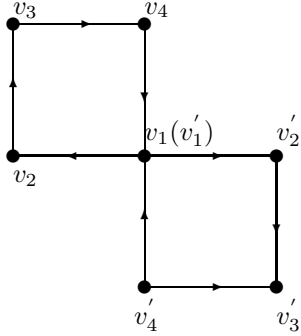
Fig.~1. A graph on $7$ vertices with an orientation.
\end{center}

\vskip 2mm

Now we give an observation on Problem 2, i.e.,  affirmative answer
when $G$ is a connected bipartite graph.

\begin{Theorem}
Let $G$ be a connected  bipartite graph. If $\rho(S)$ is the same
for all skew adjacency matrices $S$ of $G$, then $G$ is a tree.
\end{Theorem}

\begin{Proof}
Let
\[
A=\left[
\begin{array}{ccc}
 O & B  \\
 B^{\top} & O
\end{array}
\right]
\]
and
\[
\overline{S}=\left[
\begin{array}{ccc}
 O & B  \\
 -B^{\top} & O
\end{array}
\right]
\]
be an (ordinary) adjacency and a skew adjacency matrix of $G$.
 Let $S$ be a
skew adjacency matrix of $G$. Then $\rho(S)=\rho(\overline{S})$. By Lemma
\ref{lm2.5}, there is a $\{-1, 1\}$-diagonal matrix $D$
 such that $S=D\overline{S}D^{-1}$, i.e.,  $S$ is
similar to $\overline{S}$, which implies that all skew adjacency matrices
of $G$ are cospectral. Thus by Lemma \ref{lm2.2}, $G$ is an
odd-cycle graph. Since $G$ is connected and bipartite, $G$ is a
tree.
\end{Proof}

%


Let $G$ be a connected bipartite  graph on $n$ vertices. Let
$K_{r,s}$ be the complete bipartite graph with $r$ and $s$ vertices
in its two partite sets, respectively. Note that $\rho(G)<\rho(G+e)$
for an edge of the complement of $G$ (following from the
Perron-Frobenius theorem) and that
 $\rho(K_{r,s})=\sqrt{rs}$. By Lemma \ref{lm2.5}, $\rho_s(G)=\rho(G)\le \sqrt{\lfloor\frac
{n}{2}\rfloor \lceil\frac {n}{2}\rceil}$ with equality if and only
if $G=K_{\lfloor\frac {n}{2}\rfloor, \lceil\frac {n}{2}\rceil}$, cf.
\cite[Example 6.1]{Ca}.

Let $G$ be a connected  graph on $n$ vertices. Let $P_n$ be the path
on $n$ vertices. By \cite[Lemma 6.4]{Ca}, $\rho_s(G)> \rho_s(G-e)$
for an edge $e$ of $G$. Let $T$ be a spanning tree of $G$. Then by
Lemma \ref{lm2.5} and a result of Collatz and Sinogowitz \cite{CS},
$\rho_s(G)\ge \rho_s(T)=\rho (T)\ge \rho(P_n)$ with equality if and
only if $G=P_n$, cf. \cite[Example 6.3]{Ca}.

\bigskip
\noindent
{\bf Acknowledgment.} This work was supported by the Guangdong
Provincial Natural Science Foundation of China (no.~S2011010005539). The authors thank the referee for  constructive comments.

\end{document}